\documentclass[12pt,reqno]{amsart}
\usepackage[dvips]{graphics}
\usepackage{amssymb}
\usepackage{dsfont}
\usepackage{hyperref}

\numberwithin{equation}{section}

\newtheorem{thm}{Theorem}[section]
\newtheorem*{thm*}{Theorem}
\newtheorem*{thmmain*}{MAIN THEOREM}
\newtheorem{lem}{Lemma}[section]

\newtheorem{prop}[thm]{Proposition}
\newtheorem*{prop*}{Proposition}

\theoremstyle{definition}
\newtheorem{defn}{Definition}[section]
\theoremstyle{remark}
\newtheorem{rem}{Remark}[section]
\newtheorem{ex}[rem]{Example}

\newcommand{\tref}[1]{Theorem~\ref{#1}}
\newcommand{\cref}[1]{Corollary~\ref{#1}}
\newcommand{\pref}[1]{Proposition~\ref{#1}}

\newcommand{\lref}[1]{Lemma~\ref{#1}}
\newcommand{\exref}[1]{Example~\ref{#1}}

\def\diam{\mathop{\text{diam}}}
\def\rad{\mathop{\text{rad}}}

\def\dim{\mathop{\text{dim}}}

\begin{document}
\title{Centers of convex subsets of buildings}
\author{Andreas Balser and Alexander Lytchak}

\begin{abstract}
We prove that two dimensional convex subsets of spheri\-cal buildings are either
buildings or have a center. 
\end{abstract}

\maketitle
\renewcommand{\theequation}{\arabic{section}.\arabic{equation}}
\pagenumbering{arabic}

\section{Introduction}
Recently, convex cores of isometric group actions on symmetric spaces were
studied in \cite{Groups}. The results of Kleiner and Leeb imply the
following: 
Assume that the group $\Gamma$ acts by isometries on a symmetric
space~$H$ in a non-elementary way (i.e. without fixed point at
infinity). Then the 
convex core of $\Gamma$ is a direct product of  symmetric and 
Gromov-hyperbolic spaces. In particular the boundary at infinity
of the convex core must be  a spherical building with respect to the
 Tits metric.
 
This gives rise to the following more abstract question. 
Let $\Gamma$ be a group acting by isometries on a spherical building $G$.
Let $X$ be a $\pi$-convex subset of $G$ invariant under $\Gamma$. Is it
true that $\Gamma$ must have  a fixed point if $X$ is not a building?
In the case of small dimensions (of the building or at least the
subset) we answer this question in the affirmative. 

\begin{thm}\label{mainTheorem}
 Let $G$ be a spherical building, $X\subset G$ a convex subset. 
If $\dim (X) \leq 2$ then either $X$ is itself a building or
the isometry group of $X$ has a fixed point in $X$.
\end{thm}

In the proof we use some  observations about general $CAT(1)$ spaces,
 which we consider to be of 
independent interest. In important special situations they were
obtained in \cite{parabolic} by different methods. 

The next result is shown in \cite[Thm.\ B]{lang} for 
$r<\frac \pi 2$ and in \cite[Thm.\ 1.7]{parabolic} for $\diam X=\frac \pi 2$.

\begin{prop} \label{radiusDiameter}
Let $X$ be a finite-dimensional CAT(1) space. If the radius $r$ of
$X$ is smaller than $\pi$, then for the diameter  of $X$
holds $\diam (X) >\rad (X)$.
\end{prop}

From this, an easy limiting argument gives the following uniform result,
which follows from  \cite[Thm.\ B]{lang} in the case $r< \frac \pi 2$:

\begin{thm} \label{peppup}
For each natural number $n$ and each $r<\pi $ there is some 
$\varepsilon = \varepsilon (r,n) >0$, such that for each $CAT(1)$ space
$X$ of dimension at most $n$, of radius at most $r$ with at least two points,
 we have $\frac {\diam (X)} {\rad (X)} \geq 1 + \varepsilon$.  
\end{thm}

\begin{rem}
For $\frac \pi 2 \leq r < \pi$, this result can not be 
achieved by the methods of \cite{lang}. It
is not clear to us what the optimal value of $\varepsilon$ should be and 
what the spaces with smallest possible ratio 
$\frac {\diam (X)} {\rad (X)}$ look like. 
\end{rem}

In \cite[Thm.\ B]{lang} it is shown that a CAT(1) space with radius
$<\frac \pi 2$ 
has a unique circumcenter, that is therefore fixed by the whole
isometry group.  In the compact case the next result is proved in 
\cite[Prop.\ 5.7]{parabolic}.

\begin{prop} \label{existence}
Let $X$ be a  CAT(1) space of  finite dimension and of radius
$r \leq \frac{\pi}{2}$. Then $X$ has a circumcenter which is fixed by
every  isometry of $X$.
\end{prop}

\begin{rem}
In both propositions the assumption of  finite dimensionality is essential,
compare \cite[pg.\ 5]{parabolic} or \exref{hilbertSphereExample}.
\end{rem}

We explain the idea of the proof of \tref{mainTheorem} in the case
where~$X$ is compact. In this case we consider a smallest closed
convex  
subset~$Y$  of~$X$ invariant under the isometry group of $X$ and one can
assume that~$Y$ is not a building. Since $Y$ is minimal one can  not
cut off any small neighborhood of any point $y\in Y$, such that the remaining
subset is still convex (see Section \ref{action} for more on the notion
of non-removable points). In dimension $2$ this implies that each point
is an inner point of some geodesic. Now we consider two points $y,z\in Y$
with maximal distance in $Y$. If $d(y,z) =\pi$ then $Y$ contains a circle
and one can deduce that $\rad (Y) = \frac \pi 2$. If $d(y,z)<\pi$,
  we show that one can choose geodesics through $y$  resp. through
$z$. 
Using these, we construct a spherical quadrangle in which one of the
diagonals is longer than $d(y,z)$, a contradiction.
The last argument in the proof uses spherical quadrangles 
(not  triangles as usual) and breaks down if the dimension is 
bigger than $2$. 


\smallbreak

We thank Bernhard Leeb for posing the question.

\section{Preliminaries} 
\subsection{Notations} By $d$ we denote distances in metric spaces. 
By $B_r (x)$ we will denote the closed metric ball of radius $r$
around the point $x$.
For a point $x$ in a metric space $X$ we  set
 $\rad _x (X) := \sup _{z\in X} d(x,z)$. By the radius resp. the diameter
of $X$ we denote $\rad (X) := \inf _{x\in X} \rad _x (X)$ resp. 
$\diam (X):= \sup _{x\in X} \rad _x (X)$.

 Since we are dealing with non-proper spaces we will use the concept of 
ultraconvergence (instead of the Gromov-Hausdorff convergence)
with respect to a fixed non-principal ultrafilter 
$\omega$. We refer to \cite[pp.\ 77-80]{bridsonHaefliger} for
details, see also \cite[sect.\ 11]{lytchakSphB}. We will denote
by $\lim _{\omega} (X_i,x_i)$ the ultralimit of pointed metric spaces $X_i$
with respect to~$\omega$. By $X ^{\omega}$ we will denote the ultraproduct
of $X$, i.e. the ultralimit of the constant sequence $(X,x)$.

\subsection{CAT(1) spaces}
 A complete metric 
space is  called CAT(1) if each pair of 
points  with distance  $<\pi$ is connected by a geodesic and 
all triangles of perimeter less than $2\pi$ are not thicker than in
$S^2$. We refer to 
\cite[ch.\ II]{bridsonHaefliger} for the theory of such spaces.
In a CAT(1) space $X$ we will denote by $S_x$ the link at the point $x$.
A subset
$C$ of a CAT(1) space is convex if all points in $C$ with distance $<\pi$
are joined  by a geodesic in $C$.

By $\dim (X)$ we denote the geometric dimension of $X$ studied in 
\cite{kleinerDimension}.
 Ultralimits of CAT(1) spaces are CAT(1) and the dimension does not 
increase by this procedure (\cite[L.\ 11.1]{lytchakSphB}).

A subset $T$ of a CAT(1) space $X$ is  spherical if it can be 
isometrically embedded into some Euclidean sphere. See  \cite{part1} for more
on this. 

\subsection{Buildings and their convex subsets} \label{build}
We refer to \cite[sect.\ 3]{kleinerLeeb} for an account on spherical
buildings.  From the geometry of buildings we will only use the 
following basic property:

 Let $G$ be a building. Then for each $x\in X$
there is a number $r_x >0$ (set $r_{x}$ so that any point
$\bar{x}\in X$ with $d(\bar{x},x)<r_{x}$ lies in a Weyl chamber
containing $x$), 
such that for all $x,y\in X$ and all 
$\bar x \in B_{r_x} (x)$ and $\bar y \in B_{r_y} (y)$ the convex hull of the
four points $x,\bar x,y,\bar y$ is  a spherical (usually three-dimensional)
subset of $G$.  Remark that this property is inherited  by convex subsets of 
buildings. We refer to \cite{part1} for a more thorough study.

We will use the following  result derived in \cite{part1}:
\begin{lem}[{\cite[6.3 + 6.4]{part1}}]\label{part1Lemma}
Let $X$ be  a convex subset of some building $G$. If $X$ is not a building,
then $\rad (X) < \pi$. If $X$ contains an isometrically embedded $S^{n-1}$
and $\dim (X) \leq n$, then either~$X$ is a building or 
$\rad (X) = \frac \pi 2$.   
\end{lem}

\section{Radii of CAT(1) spaces}
\subsection{Generalities}
The following 
lemma follows from \cite[L.\ 11.7]{lytchakSphB},
but we give a direct proof here:
\begin{lem} \label{easyLemma}
Let $X$ be an $n$-dimensional CAT(1) space, and let $S \subset X$ be an 
embedded $S^n$. Then for each $x\in X$ there is an antipode $y\in S$, i.e.
a point satisfying $d(x,y) \geq \pi$. Therefore we have $\rad (X) \geq \pi$.
\end{lem}

\begin{proof}
If $\dim S=0$, the claim is clear. 
Choose an  arbitrary $y\in S$. Let $v \in S_y$ be the starting direction
 of $yx$. By induction, we can
choose an 
antipode~$w$ of $v$ in $S^{n-2} 
=S_yS\subset S_y X$. Now we obtain an antipode of~$x$ by extending
$xy$  inside $S$ in the
direction of $w$.
\end{proof}

The proof of the next lemma is a typical application of ultraproducts.
\begin{lem} \label{ultraproof}
Let $X$ be a CAT(1) space. If $r=\diam (X)= \rad (X) < \pi$, 
then for each point
$x\in X ^{\omega }$ holds $\rad (S_x (X^{\omega })) \leq \frac \pi 2$. 
\end{lem}

\begin{proof}
 Let $x=(x_i)$. Choose $z_i\in X$ with $d(z_i,x_i) \geq r - \frac 1 i$.
For $z=(z_i) \in X^{\omega}$ we 
get  $d(z,x)=r=\diam (X) =\diam (X^{\omega} )$.
From the CAT(1) property we immediately obtain that the starting direction
$v\in S_x$ of the geodesic $xz$ satisfies $d(v,w)\leq \frac \pi 2$
for each $w\in  S_x$. 
\end{proof}

\begin{proof}[Proof of \pref{radiusDiameter}]
In \cite[Thm.\ B]{kleinerDimension} it is shown that in each CAT(1) space $X$
of dimension   $n+1$,  some link $S_x$ contains an $n$-dimensional Euclidean
sphere. Using \lref{ultraproof} and \lref{easyLemma}, we deduce
\pref{radiusDiameter}.
\end{proof}

\begin{proof}[Proof of \tref{peppup}]
Assume the contrary and choose a sequence of spaces $X_i$ with
$\frac {\diam(X_i)} {\rad (X_i)} \to 1$. Since we assume that the spaces
 have  more than one point we may rescale $X_i$ (without leaving the category
of at most $n$-dimensional $CAT(1)$ spaces)
 and assume that the radius of all the spaces $X_i$ is equal to $r$. For the
ultralimit space  $X= \lim _{\omega} (X_i)$ we obtain $\rad (X) =\diam (X)=r$
in contradiction to \pref{radiusDiameter}, since the dimension of $X$ is 
bounded by $n$ too. 
\end{proof}

\begin{ex}\label{hilbertSphereExample}
 The subset $P$
of points in the Hilbertsphere with all
coordinates  non-negative satisfies $\rad (P)=\diam (P)= \frac \pi 2$. 
Moreover, since swapping coordinates is an isometry, the isometry group of $P$ 
has no fixed points, thus the assumption of finite-dimensionality is 
essential in \pref{radiusDiameter} and \pref{existence}.  
(This was also observed in  \cite[pg.\ 5]{parabolic}).
Moreover  $P$ contains points $x$ with $\rad (S_x P)=\pi$ (let~$x$ be a
point with all components positive). 
\end{ex}

\subsection{CAT(1) spaces of radius $\frac \pi 2$} 
We recall

\begin{defn}
A point $x$ in a metric space $X$ is a \emph{circumcenter} of $X$ if
$\rad_{x}(X)=\rad X$ holds.
\end{defn}

 Due to \cite[Thm.\ B]{lang} a CAT(1) space $X$ of radius $< \frac
\pi 2$ has a unique 
circumcenter. As a first step, we extend the existence result:

\begin{lem}\label{circumcentersExist}
Let $X$ be a CAT(1) space of finite dimension and
radius~$\frac{\pi}{2}$.  Then 
$X$ has a circumcenter.
\end{lem}
\begin{proof}
By definition, we find a sequence of points $x_{i}\in X$ such that
$\rad_{x_{i}}(X) \rightarrow \frac{\pi}{2}$. This sequence defines a
circumcenter $x =(x_i)$
 for the ultraproduct $X^{\omega}$, which is a space of
the same dimension as $X$ (\cite[Cor. 11.2]{lytchakSphB}). We have
$d(x,X)\leq \frac{\pi}{2}$. In  the case of equality we have 
$d(x,y)= \frac{\pi}{2}$ for all
$y\in X$, and the convex hull of $x$ and $X$ in $X^{\omega}$
 would be isometric to the spherical join 
$X*\{x \}$ (\cite[L.\ 4.1]{lytchakSphB}). This contradicts 
$\dim (X^{\omega} )= \dim (X) $.

Hence, there is a unique projection $x'$ of $x$ to $X$ 
(\cite[II.2.6.1]{bridsonHaefliger}).
 Triangle comparison shows that $x'$ is a circumcenter of $X$.
\end{proof}

\begin{proof}[Proof of \pref{existence}]
The proof can be concluded in the same manner as in
the proof of \cite[Thm.\ 1.3]{parabolic}. The set $X_1$ of all
circumcenters of $X$ is  (by the above)  
a non-empty subset of $X$, which is closed, convex and invariant under
the isometry group. Moreover by definition the diameter of $X_1$ is not
not bigger than  $\rad (X)= \frac \pi 2$.

Due to \pref{radiusDiameter}, the radius of $X_1$
is smaller than $\frac \pi 2$ and due to \cite[Thm.\ B]{lang}, the
space $X_1$ has  
a unique circumcenter, that is therefore 
fixed by the isometry group of $X$.
\end{proof}

\section{Groups acting on a CAT(1) space} \label{action}
\subsection{Removable points} For $r>0$ we will denote by $S_r (x)$ the
sphere of radius $r$ around $x$, i.e. the set of all points $y$ with 
$d(y,x) =r$.

\begin{rem}
Do not confuse the sphere $S_r (x)$ with the link $S_x$.
\end{rem}

\begin{defn}
Let $X$ be a CAT(1) space and $r>0$. A
 point $x\in X$ is called \emph{$r$-removable}, if
 the closed convex hull of the sphere $S_r (x)$ does not
contain~$x$.  It is called \emph{removable} if it is $r$-removable for all
$r>0$. A point which is not removable is called \emph{non-removable}.
\end{defn}

Remark that if $x\in X$ is $r$-removable then it is $r'$-removable for
all $r' >r$. A point $x\in X$ is removable iff it is removable inside some
ball $B_r (x)$.  The following examples should only illustrate the notion
and will not be used in the proof below.

\begin{ex}
In every geodesically complete CAT(1) space each point is non-removable.
\end{ex}

\begin{ex}
A point $x\in X$ is non-removable  if it
is inner point of some geodesic or, more general,  if
it is an inner point of some geodesic in the ultraproduct $X^{\omega }$. 
If $X$ has dimension $1$ the last condition is also necessary.
\end{ex}

\begin{ex} 
For every CAT(1) space $Z$, each point in the spherical join
$X=Z*S^1$ is a midpoint of a 
geodesic  of length $\pi$, hence it is not $r$-removable for 
$r<\frac \pi 2$.
\end{ex}
\begin{ex} 
A closed convex subset $S$ of a finite-dimensional sphere has no removable
points iff $S$ has the form $S= Z*S^1$.
\end{ex}

\begin{ex}
Let $X$ be a locally conical space (\cite[sect.\ 3.2]{part1})
 with $\dim (X) \leq 2$. A point
$x\in X$ is non-removable iff $x$ is an inner point of some geodesic 
$\gamma \subset X^{\omega}$ in the ultraproduct $X^{\omega}$ of $X$.
\end{ex}

\begin{ex} 
 Let $X$ be a closed convex subset in a Riemannian mani\-fold of
constant curvature. A point $x\in X$
is not removable iff it is the inner point of some geodesic.
An interesting question which we  conjecture to be true is whether
this statement is also true if $X$ 
is a Riemannian manifold with variable curvature.
\end{ex}

Unfortunately it is difficult to say much about the behavior of convex hulls 
and therefore of removable points under ultralimits. The following 
easy observation will be used below
\begin{lem} \label{ultrahull}
Let $(X_i,o_i)$ be CAT(1) spaces, $(X,o)=\lim _ {\omega} (X_i,o_i)$.
Consider a point $x=(x_{i})\in X$ such that for some fixed $s>0$
the point~$x_i$ is not $s$-removable  
in $X_i$.
Then for each $m\in X$ with $l= d(m,x) \leq \frac \pi 2$, there is a point
$q\in X$ with $d(q,x) =s$ and $d(m,q) \geq d(m,x)$. 
\end{lem}

\begin{proof}
We can represent $m$ as a sequence of points $m=(m_i)$, with $m_i \in X_i$,
$l_i =d(m_i,x_i) < \frac \pi 2$. We may assume that $m\not =x$ and
hence, that $l_{i}\not \rightarrow 0$. Since $x_i$ is not
$s$-removable, there must be at least one point $q_i$ with
$d(q_i,x_i)=s$, that is not contained in the 
closed convex ball $B_{l_i-\rho _i} (m_i)$, for each $\rho _i >0$. If
we choose $\rho _i \to 0$, then the point $q=(q_i) \in X$ has the 
desired properties.
\end{proof}

\subsection{Minimal invariant subsets} Non-removable points are related to
 minimal invariant subsets by the following observation:

\begin{lem}\label{nonRemovablePoints}
Let $X$ be a 
CAT(1) space, $x\in X$, $\Gamma$ a group of isometries of $X$.
Assume that there is no proper closed, convex, and $\Gamma$-invariant
subset of $X$ 
containing~$x$. 
Then each removable point of $X$ is contained in
the closure of the orbit $\overline{\Gamma x}$. More precisely if
$d(z,\Gamma x) =\varepsilon$, then $z$ is not $r$-removable
 for $r<\varepsilon$.
\end{lem}

\begin{proof}
Assume the contrary and 
 consider the closed convex convex hull  $C$ of $X\setminus B_r (z)$.
 Then $C$ does not contain $z$ but  
$x\in g C$ for $g\in \Gamma$. Considering $\cap _{g \in \Gamma } gC$ we
get a contradiction to the minimality of $X$.
\end{proof}

\begin{rem}
If $X$ is a minimal closed convex subset invariant under $\Gamma$, then
either $X$ has no removable points or $\Gamma$ has a dense orbit. If $X$
is finite dimensional and not discrete one can show that in the latter
case each point is  
non-removable as well.  
\end{rem}

\subsection{Centers} We will stick to the following
\begin{defn}
We say that a CAT(1) space $X$ has a \emph{center} if 
its group of isometries
$\Gamma ^X =Iso (X)$ has a fixed point.
\end{defn}

 A group operating isometrically on a CAT(1) space $X$ has a fixed point 
iff it has an orbit of diameter $< \frac \pi 2$. This implies

\begin{lem} \label{CFU}
 Let $X_i$ be a sequence of CAT(1) spaces. If $X_i$ has no center,
then $X= \lim _{\omega } X_i$ has no center too. 
\end{lem}

\section{The Proof of the Main Theorem} \label{proofSection}
For  the proof of \tref{mainTheorem} the following lemma is crucial:

\begin{lem}\label{diamLemma}
Let $Y$ be a convex subset of a  building, $\dim (Y) \leq 2$. Let
$y,z \in Y$ be inner points of some geodesics in $Y$. If $d(y,z)=\diam (Y)$
holds, then $d(y,z)=\pi$ and $Y$ contains some $S^1$.
\end{lem} 

\begin{proof}
Let $y$ resp. $z$ be inner points of geodesics $p_1p_2$ resp. $q_1q_2$.
We may assume that $p_i \in B_{r_y} (y)$ resp. $q_i \in B_{r_z} (z)$ for 
 $i=1,2$, where $r_y$ resp. $r_z$ are the conicality radii at $y$ resp. at
 $z$ (see Subsection \ref{build}).
 If $d(y,z)=\pi$, the convex hull of $q_1q_2$ and $y$ 
is  a circle  $S^1$. 

Thus assume that $d(y,z)<\pi$. For the geodesic $\eta =yz$
we deduce from the assumption $d(y,z)= \diam (Y)$, that 
$\eta$ meets   $yp_i$ and $zq_j$  orthogonally.
From the geometry of buildings we derive (cf.\ Subsection~\ref{build} again)
that the convex hulls of $yp_i$
and $zq_j$ are spherical for $i,j=1,2$. Since $\dim (Y) \leq 2$ these
spherical hulls are 2-dimensional spherical quadrangles.

 Now in the two dimensional sphere the following statement holds:
For $i=1,2$ let  $\bar \gamma _i$ be geodesics in $S^2$ starting in $x_i$
with $d(x_1,x_2) <\pi$.
Assume that $\bar \eta =x_1x_2$  meets $\bar \gamma _1$ and $\bar \gamma _2$
orthogonally. If $\bar \gamma _1$ and $\bar \gamma _2$ are on different sides
of $\eta$, i.e.  
if $\angle _{x_1} (\bar \gamma _1 (\varepsilon ), \bar \gamma _2 (\varepsilon))
>\frac \pi 2$, 
then $d(\bar \gamma _1 (\varepsilon ), \bar \gamma _2 (\varepsilon ) )
>d(x_1,x_2)$
for small $\varepsilon >0$.

 Therefore we get a contradiction to 
$d(y,z)= \diam (Y)$ as soon as we can verify that 
$\angle _y (p_i,q_j)  >\frac \pi 2$ for some $i$ and $j$.  However
from the one-dimensionality of $S_y Y$ we directly obtain
that $\angle _y (p_1,q_j) >\frac \pi 2$ for $j=1$ or $j=2$.
\end{proof}

Now we are going to prove \tref{mainTheorem}
 along the line explained in the introduction.
The non-compactness causes additional difficulties.

\begin{proof}[Proof of \tref{mainTheorem}] Let $G$ be a spherical
building modelled on the Coxeter group~$W$. Let $X \subset G$ be a
convex subset of dimension at most~$2$ and without center. Assume
finally that $X$ is not a building. Due to \lref{part1Lemma} we have
$\bar r := \rad (X) <\pi$.

Let  $C _W ^{\bar{r}}$ be the class of all CAT(1) spaces that have dimension 
at most $2$, no center, whose radius is at most $\bar r$ and that admit
an iso\-metric embedding into some spherical building modelled on the Coxeter
group~$W$.
 Since an ultralimit of buildings modelled on $W$ is again a building 
modelled on $W$, the class  $C _W ^{\bar{r}}$ is closed under
ultralimits (by \lref{CFU}).  
 Therefore we can find a space $Z$ in $C _W ^{\bar{r}}$,
whose radius 
$r\leq \bar r$ is the smallest possible for spaces in  $C _W
^{\bar{r}}$ (let $r_{i}\rightarrow r$ be radii of elements of
$C_{W}^{\bar{r}}$, and let $Z$ be the ultralimit of the corresponding spaces).

Replacing $Z$ by its ultraproduct $Z ^{\omega}$ we may assume that
there is a point $x\in Z$ with $\rad _x (Z) = \rad (Z)$. Let $Y$ be 
the smallest closed convex subset of $Z$ containing $x$ and invariant under
the isometry group $\Gamma ^Z$ of $Z$. 
Clearly, a center for $Y$ would also be a center for $Z$, thus $Y$ is in
$C^{r}_{W}$ too. Moreover we have $\rad Y= \rad _x (Y) =r$ by the 
minimality of $r$, and $r>\frac \pi 2$ by Prop. \ref{existence} since $Y$ has no center.

From  \pref{radiusDiameter} we derive  $\diam Y -\rad _x (Y)
=2\varepsilon  >0$.  Therefore for 
points  $y_i,z_i \in Y$ with $d(y_i,z_i) \to 
\diam (Y)$ we can apply Lemma \ref{nonRemovablePoints} and see that
for all big $i$ the points 
$y_i$ and $z_i$ are not $\varepsilon $-removable.

 Consider $y=(y_i), z=(z_i) \in Y ^{\omega}$. 
We have $d(y,z)=\diam Y^{\omega}=\diam Y$. Assume that $y$ and $z$
are inner points of some geodesics in~$ Y^{\omega}$. Then from
\lref{diamLemma} we derive that $Y^{\omega}$ contains a circle, from 
\lref{part1Lemma}
we see that $\rad (Y^{\omega }) =\frac \pi 2$ and by  \pref{existence}
 the space  $Y^{\omega}$ must have a center - a contradiction.

By symmetry it is enough to prove that $y$ is an inner point of a geodesic.
Decreasing $\varepsilon$  we may assume that it
is smaller than the conicality radius $r_y$ of $y$ in $Y^{\omega}$.

Denote by $\eta$ the geodesic between $y$ and $z$ and let $v\in S_y$
be its starting direction. From   $d(y,z)=\diam Y^{\omega}$ we derive 
$\rad _v (S_y) \leq \frac \pi 2 $.

Set $m=\eta (\frac \pi 2 )$. By \lref{ultrahull} there is a point $q$ with 
$d(q,m) \geq \frac \pi 2$ and $d(y,m)= \varepsilon$. Since the triangle $ymq$
is spherical,  for the starting direction $w$ of $yq$ we obtain
$d(v,w)= \frac \pi 2$.

Now we consider the point $z_{\delta}$ on the geodesic $qz$ with
$d(z,z_{\delta})= \delta \to 0$. As above we derive from \lref{ultrahull}
that there is some $q_{\delta} \in Y^{\omega}$ such that
$d(q_{\delta} ,y)=\varepsilon$ and such that for the starting directions
$v_{\delta} $ of $yz_{\delta}$ resp. $w_{\delta} $ of $yq_{\delta}$ the 
inequality $d(v_{\delta},w_{\delta}) \geq \frac \pi 2$ holds.
 
Since the triangle $yzq$ is spherical we know $d(w,v_{\delta}) +
d(v_{\delta}, v)= d(v,w) = \frac{\pi}{2}$. Therefore from $\rad _v
(S_y) =\frac \pi 2$ and 
the fact that $S_y Y^{\omega}$ is one-dimensional we deduce that 
$d(w,w_{\delta} )\to \pi$. Now using a diagonal argument in $Y$ 
(or going to another ultraproduct $(Y^{\omega})^{\omega}$) we obtain
a point $q_0$ such that $y$ is an inner point of the geodesic $qq_0$.

This finishes the proof of \tref{mainTheorem}.
\end{proof}

\bibliographystyle{alpha}
\bibliography{fix}

\end{document}